\documentclass[10pt,letterpaper]{article}
\usepackage{theorem,enumerate}
\usepackage{amsmath,latexsym,amssymb,amsfonts}
\usepackage{eucal}
\usepackage{color}
\usepackage{comment}
\theorembodyfont{\normalfont\slshape}
\newtheorem{thm}{Theorem}[section]
\newtheorem{cor}[thm]{Corollary}
\newtheorem{prop}[thm]{Proposition}
\newtheorem{lem}[thm]{Lemma}

\newtheorem{claim}{Claim}[section]

\newtheorem{Thm}{Theorem}

\newtheorem{remark}{Remark}

\newcommand{\proof}{\medbreak\noindent\textit{Proof.}\quad}

\newcommand{\qed}{{$\quad\square$\vs{3.6}}}

\newcommand{\vs}[1]{\vspace*{#1 mm}}

\numberwithin{equation}{section}

\def\MM{{ \mathcal{M}}}

\def\SS{{ \mathcal{S}}}

\bfseries\normalfont

\title{Small domination-type invariants in random graphs}
\author{
Michitaka Furuya\footnote{College of Liberal Arts and Science, Kitasato University, 1-15-1 Kitasato, Minami-ku, Sagamihara, Kanagawa 252-0373, Japan. \texttt{e-mail:michitaka.furuya@gmail.com}},
Tamae Kawasaki\footnote{Department of Applied Mathematics, Tokyo University of Science, 1-3 Kagurazaka, Shinjuku-ku, Tokyo 162-8601, Japan. \texttt{e-mail:tm.kawasaki@rs.tus.ac.jp}}
}

\date{}

\begin{document}

\maketitle

\begin{abstract}
For $c\in \mathbb{R}^{+}\cup \{\infty \}$ and a graph $G$, a function $f:V(G)\rightarrow \{0,1,c\}$ is called a $c$-self dominating function of $G$ if for every vertex $u\in V(G)$, $f(u)\geq c$ or $\max\{f(v):v\in N_{G}(u)\}\geq 1$ where $N_{G}(u)$ is the neighborhood of $u$ in $G$.
The minimum weight $w(f)=\sum _{u\in V(G)}f(u)$ of a $c$-self dominating function $f$ of $G$ is called the $c$-self domination number of $G$.
The $c$-self domination concept is a common generalization of three domination-type invariants; (original) domination, total domination and Roman domination.
In this paper, we study a behavior of the $c$-self domination number in random graphs for small $c$.
\end{abstract}



%%%%%%%%%%%%%%%%%%%%%%%%%%%%%%%%%%%%%%%%%%%%%%%%%%%%%%%%%%%%%%%%%%%%%%%%%%%%%%%%%%%%%%%%%%%%%%%%%%%%%%%%%%%%%%%%%%%%%%%%
%%%%%%%%%%%%%%%%%%%%%%%%%%%%%%%%%%%%%%%%%%%%%%%%%%%%%%%%%%%%%%%%%%%%%%%%%%%%%%%%%%%%%%%%%%%%%%%%%%%%%%%%%%%%%%%%%%%%%%%%
%%%%%%%%%%%%%%%%%%%%%%%%%%%%%%%%%%%%%%%%%%%%%%%%%%%%%%%%%%%%%%%%%%%%%%%%%%%%%%%%%%%%%%%%%%%%%%%%%%%%%%%%%%%%%%%%%%%%%%%%
\section{Introduction}\label{sec1}
%%%%%%%%%%%%%%%%%%%%%%%%%%%%%%%%%%%%%%%%%%%%%%%%%%%%%%%%%%%%%%%%%%%%%%%%%%%%%%%%%%%%%%%%%%%%%%%%%%%%%%%%%%%%%%%%%%%%%%%%
%%%%%%%%%%%%%%%%%%%%%%%%%%%%%%%%%%%%%%%%%%%%%%%%%%%%%%%%%%%%%%%%%%%%%%%%%%%%%%%%%%%%%%%%%%%%%%%%%%%%%%%%%%%%%%%%%%%%%%%%
%%%%%%%%%%%%%%%%%%%%%%%%%%%%%%%%%%%%%%%%%%%%%%%%%%%%%%%%%%%%%%%%%%%%%%%%%%%%%%%%%%%%%%%%%%%%%%%%%%%%%%%%%%%%%%%%%%%%%%%%

Throughout this paper, we let $\mathbb{R}^{+}$ and $\mathbb{Z}^{+}$ denote the set of positive numbers and the set of positive integers, respectively.
Let $G$ be a graph.
Let $V(G)$ and $E(G)$ denote the vertex set and the edge set of $G$, respectively.
For a vertex $u\in V(G)$, we let $N_{G}(u)$ denote the {\it neighborhood} of $u$ in $G$; thus $N_{G}(u)=\{v\in V(G): uv\in E(G)\}$.
A set $S\subseteq V(G)$ is a {\it dominating set} (resp. a {\it total dominating set}) of $G$ if each vertex in $V(G)\setminus S$ (resp. each vertex in $V(G)$) is adjacent to a vertex in $S$.
The minimum size of a dominating set (resp. a total dominating set) of $G$, denoted by $\gamma (G)$ (resp. $\gamma _{t}(G)$), is called the {\it domination number} (resp. the {\it total domination number}) of $G$.
Since a graph $G$ with isolated vertices has no total dominating set, the total domination number has been typically defined for only graphs without isolated vertices.
However, in this paper, we define $\gamma _{t}(G)$ as $\gamma _{t}(G)=\infty $ if $G$ has an isolated vertex for convenience.
Domination and total domination are important invariants in graph theory because they have many applications for mathematical problems and real problems (see \cite{HHS1,HHS2,HY}).

%Many generalizations of the domination number have been known, for example, $k$-domination, distance domination, rainbow domination, etc.
The first author~\cite{F1} recently defined a new domination-type concept as follows:
Let $G$ be a graph.
For a function $f:V(G)\rightarrow \mathbb{R}^{+}\cup \{0,\infty \}$, the {\it weight} $w(f)$ of $f$ is defined by $w(f)=\sum _{u\in V(G)}f(u)$.
Let $c\in \mathbb{R}^{+}\cup \{\infty \}$.
A function $f:V(G)\rightarrow \mathbb{R}^{+}\cup \{0,\infty \}$ is a {\it $c$-self dominating function} (or {\it $c$-SDF}) of $G$ if for each $u\in V(G)$, $f(u)\geq c$ or $\max\{f(v):v\in N_{G}(u)\}\geq 1$.
Then the following proposition holds.

\begin{prop}[Furuya~\cite{F1}]%%%%%%%%%%%%%%%%%%%%%%%%%%%%%%%%%%%%%%%%%%%%%%%%%%%%%%%%%%%%%%%%%%%%%%%%%%%%%%%%%%%%%%%%%%
\label{prop1.1}
Let $c\in \mathbb{R}^{+}\cup \{\infty \}$, and let $G$ be a graph.
If $f$ is a $c$-SDF of $G$, then there exists a $c$-SDF $g$ of $G$ such that $w(g)\leq w(f)$ and $g(u)\in \{0,1,c\}$ for all $u\in V(G)$.
\end{prop}
%%%%%%%%%%%%%%%%%%%%%%%%%%%%%%%%%%%%%%%%%%%%%%%%%%%%%%%%%%%%%%%%%%%%%%%%%%%%%%%%%%%%%%%%%%%%%%%%%%%%%%%%%%%%%%%%%%%%%%%%

It follows from Proposition~\ref{prop1.1} that the minimum weight of a $c$-SDF of $G$ is well-defined.
The minimum weight of a $c$-SDF of $G$, denoted by $\gamma ^{c}(G)$, is called the {\it $c$-self domination number} of $G$.
%A $c$-SDF $f$ of $G$ with $w(f)=\gamma ^{c}(G)$ is called a {\it $\gamma ^{c}$-function} of $G$.
Note that $\gamma ^{1}(G)=\gamma (G)$ and $\gamma ^{\infty }(G)=\gamma _{t}(G)$ for all graphs $G$ (see \cite{F1}).
Furthermore, the $\frac{1}{2}$-self domination number is equal to the half of the Roman domination number defined in Subsection~\ref{sec1.1}.
Thus self domination concept is a common generalization of three well-studied invariants.

In this paper, our main aim is to analyze a behavior of the $c$-self domination number in Erd\H{o}s-R\'{e}nyi model random graphs $G(n,p)$ on $[n]:=\{1,2,\ldots ,n\}$.
For $p\in (0,1)$ and $n\in \mathbb{Z}^{+}\setminus \{1\}$, let $a_{p}(n)=\log _{1/(1-p)}\frac{n}{\log _{1/(1-p)}n\ln n}$.
Then the following are known.

\begin{Thm}[Wieland and Godbole~\cite{WG}]%%%%%%%%%%%%%%%%%%%%%%%%%%%%%%%%%%%%%%%%%%%%%%%%%%%%%%%%%%%%%%%%%%%%%%%%%%%%%%
\label{ThmA}
For $p\in (0,1)$, $\gamma (G(n,p))\in \{\lfloor a_{p}(n) \rfloor +1,\lfloor a_{p}(n) \rfloor +2\}$ with probability that tend to $1$ as $n\rightarrow \infty $.
\end{Thm}
%%%%%%%%%%%%%%%%%%%%%%%%%%%%%%%%%%%%%%%%%%%%%%%%%%%%%%%%%%%%%%%%%%%%%%%%%%%%%%%%%%%%%%%%%%%%%%%%%%%%%%%%%%%%%%%%%%%%%%%%

\begin{Thm}[Bonato and Wang~\cite{BW}]%%%%%%%%%%%%%%%%%%%%%%%%%%%%%%%%%%%%%%%%%%%%%%%%%%%%%%%%%%%%%%%%%%%%%%%%%%%%%%%%%%
\label{ThmB}
For $p\in (0,1)$, $\gamma _{t}(G(n,p))\in \{\lfloor a_{p}(n) \rfloor +1,\lfloor a_{p}(n) \rfloor +2\}$ with probability that tend to $1$ as $n\rightarrow \infty $.
\end{Thm}
%%%%%%%%%%%%%%%%%%%%%%%%%%%%%%%%%%%%%%%%%%%%%%%%%%%%%%%%%%%%%%%%%%%%%%%%%%%%%%%%%%%%%%%%%%%%%%%%%%%%%%%%%%%%%%%%%%%%%%%%

\begin{remark}%%%%%%%%%%%%%%%%%%%%%%%%%%%%%%%%%%%%%%%%%%%%%%%%%%%%%%%%%%%%%%%%%%%%%%%%%%%%%%%%%%%%%%%%%%%%%%%%%%%%%%%%%%
\label{remark1}
Recall that our definition of total domination is not traditional because we define $\gamma _{t}(G)=\infty $ for graphs $G$ with an isolated vertex.
Thus, strictly speaking, total domination in Theorem~\ref{ThmB} is different from one in this paper.
However, Bonato and Wang~\cite{BW} indeed proved that $G(n,p)$ has a total dominating set having size $\lfloor a_{p}(n) \rfloor +2$ with probability that tend to $1$ as $n\rightarrow \infty $.
Furthermore, since $\gamma (G)\leq \gamma _{t}(G)$ for all graphs $G$, it follows from Theorem~\ref{ThmA} that $G(n,p)$ has no total dominating set having the size $\lfloor a_{p}(n) \rfloor $ with probability that tend to $1$ as $n\rightarrow \infty $.
Hence Theorem~\ref{ThmB} holds under our definition.
\end{remark}
%%%%%%%%%%%%%%%%%%%%%%%%%%%%%%%%%%%%%%%%%%%%%%%%%%%%%%%%%%%%%%%%%%%%%%%%%%%%%%%%%%%%%%%%%%%%%%%%%%%%%%%%%%%%%%%%%%%%%%%%

By the definition of self domination, if $c,c'\in \mathbb{R}^{+}\cup \{\infty \}$ satisfy $c\leq c'$, then $\gamma ^{c}(G)\leq \gamma ^{c'}(G)$ for all graphs $G$.
Here we note that for $c~\in (1,\infty )$, the value $\gamma ^{c}(G)$ may be a non-integer if $c$ is a non-integer.
Thus the following result is obtained as a corollary of Theorems~\ref{ThmA} and \ref{ThmB}.

\begin{cor}%%%%%%%%%%%%%%%%%%%%%%%%%%%%%%%%%%%%%%%%%%%%%%%%%%%%%%%%%%%%%%%%%%%%%%%%%%%%%%%%%%%%%%%%%%%%%%%%%%%%%%%%%%%%%
\label{corAB}
For $c\in [1,\infty )$ and $p\in (0,1)$, $\gamma ^{c}(G(n,p))\in [\lfloor a_{p}(n) \rfloor +1,\lfloor a_{p}(n) \rfloor +2]$ with probability that tend to $1$ as $n\rightarrow \infty $.
\end{cor}
%%%%%%%%%%%%%%%%%%%%%%%%%%%%%%%%%%%%%%%%%%%%%%%%%%%%%%%%%%%%%%%%%%%%%%%%%%%%%%%%%%%%%%%%%%%%%%%%%%%%%%%%%%%%%%%%%%%%%%%%

In this paper, we focus on $c$-self domination in the remaining case, that is, the case where $c\in (0,1)$.
To state our main result, we extend the floor $\lfloor * \rfloor $.
For $t\in \mathbb{Z}^{+}$ and $a\in \mathbb{R}$, let $\lfloor a \rfloor _{t}$ be the largest number in $\{m_{1}+\frac{m_{2}}{t}:m_{1},m_{2}\in \mathbb{Z},~m_{1}+\frac{m_{2}}{t}\leq a\}$.
Recall that $a_{p}(n)=\log _{1/(1-p)}\frac{n}{\log _{1/(1-p)}n\ln n}$.
For $p\in (0,1)$, $t\in \mathbb{Z}^{+}$ and $n\in \mathbb{Z}^{+}\setminus \{1\}$, let $b_{p,t}(n)=\lfloor \lfloor a_{p}(n) \rfloor _{t}+\frac{1}{t} \rfloor +1$.
Note that if $\lfloor a_{p}(n) \rfloor _{t}+\frac{1}{t}$ is a non-integer, then $b_{p,t}(n)$ is the smallest integer more than $a_{p}(n)$; if $\lfloor a_{p}(n) \rfloor _{t}+\frac{1}{t}$ is an integer, then $b_{p,t}(n)$ is the second smallest integer more than $a_{p}(n)$.
Our main result is the following:

\begin{thm}%%%%%%%%%%%%%%%%%%%%%%%%%%%%%%%%%%%%%%%%%%%%%%%%%%%%%%%%%%%%%%%%%%%%%%%%%%%%%%%%%%%%%%%%%%%%%%%%%%%%%%%%%%%%%
\label{thm-main}
Let $s$ and $t$ be integers with $2\leq s\leq t-1$.
Then for $p\in (0,1)$,
$$
\gamma ^{\frac{t}{s}}(G(n,p))\in \left[\lfloor a_{p}(n) \rfloor _{t}+\frac{1}{t},~b_{p,t}(n)\right]\setminus \left\{b_{p,t}(n)-\frac{i}{t}:t-s+1\leq i\leq t-1\right\}
$$
with probability that tend to $1$ as $n\rightarrow \infty $.
\end{thm}
%%%%%%%%%%%%%%%%%%%%%%%%%%%%%%%%%%%%%%%%%%%%%%%%%%%%%%%%%%%%%%%%%%%%%%%%%%%%%%%%%%%%%%%%%%%%%%%%%%%%%%%%%%%%%%%%%%%%%%%%

Modeling on existing researches, we find a random variable corresponding to $c$-SDFs and calculate its expected value in Section~\ref{sec3}.
Then we will obtain a weaker result than Theorem~\ref{thm-main}: 
$$
{\rm Pr}\left(\gamma ^{\frac{t}{s}}(G(n,p))\in \left[\lfloor a_{p}(n) \rfloor _{t}+\frac{1}{t},~b_{p,t}(n)\right]\right)\rightarrow 1~~(n\rightarrow \infty )
$$
(see Theorem~\ref{thm3.0}).
The highlight of this paper is Section~\ref{sec4}.
While many known results for domination-type invariants in random graphs are completed by just calculating of a random variable, we can refine the above weak result to Theorem~\ref{thm-main} using additional graph-theoretic approach.
Note that $b_{p,t}(n)\leq \lfloor a_{p}(n) \rfloor _{t}+\frac{t+1}{t}$ and $\gamma ^{\frac{s}{t}}(G)\in \{m_{1}+\frac{m_{2}}{t}:m_{1},m_{2}\in \mathbb{Z}^{+}\cup \{0\}\}$ for all graphs $G$.
Thus Theorem~\ref{thm3.0} claims that $\gamma ^{\frac{s}{t}}(G(n,p))$ takes at most $t+1$ values with high probability, and Theorem~\ref{thm-main} improves ``at most $t+1$'' to ``at most $t-s+2$''.
In Subsection~\ref{sec1.1}, we focus on the Roman domination number an its related topic.

\begin{remark}%%%%%%%%%%%%%%%%%%%%%%%%%%%%%%%%%%%%%%%%%%%%%%%%%%%%%%%%%%%%%%%%%%%%%%%%%%%%%%%%%%%%%%%%%%%%%%%%%%%%%%%%%%
\label{remark2}
Using similar strategy in Sections~\ref{sec3} and \ref{sec4}, we can estimate $\gamma ^{c}(G(n,p))$ even if $c\in (0,1)$ is irrational number.
However, it seems to be difficult to describe an optimal formula.
On the other hand, we can give the following rough formula (by Theorem~\ref{thm3.0}):
Let $c\in (0,1)$ be an irrational number.
Then for $p\in (0,1)$ and $\varepsilon \in \mathbb{R}^{+}$, ${\rm Pr}(\gamma ^{c}(G(n,p))\in (a_{p}(n),a_{p}(n)+1+\varepsilon ])\rightarrow 1~~(n\rightarrow \infty )$.
\end{remark}
%%%%%%%%%%%%%%%%%%%%%%%%%%%%%%%%%%%%%%%%%%%%%%%%%%%%%%%%%%%%%%%%%%%%%%%%%%%%%%%%%%%%%%%%%%%%%%%%%%%%%%%%%%%%%%%%%%%%%%%%

%%%%%%%%%%%%%%%%%%%%%%%%%%%%%%%%%%%%%%%%%%%%%%%%%%%%%%%%%%%%%%%%%%%%%%%%%%%%%%%%%%%%%%%%%%%%%%%%%%%%%%%%%%%%%%%%%%%%%%%%
%%%%%%%%%%%%%%%%%%%%%%%%%%%%%%%%%%%%%%%%%%%%%%%%%%%%%%%%%%%%%%%%%%%%%%%%%%%%%%%%%%%%%%%%%%%%%%%%%%%%%%%%%%%%%%%%%%%%%%%%
\subsection{Roman domination and differential}\label{sec1.1}
%%%%%%%%%%%%%%%%%%%%%%%%%%%%%%%%%%%%%%%%%%%%%%%%%%%%%%%%%%%%%%%%%%%%%%%%%%%%%%%%%%%%%%%%%%%%%%%%%%%%%%%%%%%%%%%%%%%%%%%%
%%%%%%%%%%%%%%%%%%%%%%%%%%%%%%%%%%%%%%%%%%%%%%%%%%%%%%%%%%%%%%%%%%%%%%%%%%%%%%%%%%%%%%%%%%%%%%%%%%%%%%%%%%%%%%%%%%%%%%%%

A function $f:V (G)\rightarrow \{0,1,2\}$ is a {\it Roman dominating function} of $G$ if each vertex $u\in V(G)$ with $f(u)=0$ is adjacent to a vertex $v\in V(G)$ with $f(v)=2$.
The minimum weight of a Roman dominating function of $G$, denoted by $\gamma _{R}(G)$, is called the {\it Roman domination number} of $G$.
Roman domination was introduced by Stewart~\cite{S}, and was studied by Cockayne et al.~\cite{CDHH} in earnest.
Since $\gamma _{R}(G)=2\gamma ^{\frac{1}{2}}(G)$ for all graphs $G$, we obtain the following result as a corollary of Theorem~\ref{thm-main}.

\begin{cor}%%%%%%%%%%%%%%%%%%%%%%%%%%%%%%%%%%%%%%%%%%%%%%%%%%%%%%%%%%%%%%%%%%%%%%%%%%%%%%%%%%%%%%%%%%%%%%%%%%%%%%%%%%%%%
\label{cor-Roman}
For $p\in (0,1)$, $\gamma _{R}(G(n,p))\in \left\{2\lfloor a_{p}(n) \rfloor _{2}+i:1\leq i\leq 3\right\}$ with probability that tend to $1$ as $n\rightarrow \infty $.
\end{cor}
%%%%%%%%%%%%%%%%%%%%%%%%%%%%%%%%%%%%%%%%%%%%%%%%%%%%%%%%%%%%%%%%%%%%%%%%%%%%%%%%%%%%%%%%%%%%%%%%%%%%%%%%%%%%%%%%%%%%%%%%

Roman domination is closely related to another important invariant.
The {\it differential} of a graph $G$, denoted by $\partial (G)$, is defined as $\partial (G)=\max\{|(\bigcup _{u\in X}N_{G}(u))-X|-|X|:X\subseteq V(G)\}$.
The differential has been widely studied because it was motivated from information diffusion in social networks.
Recently, Bermudo et al.~\cite{BFS} proved a very useful result that every graph $G$ satisfies $\gamma _{R}(G)+\partial (G)=|V(G)|$.
Thus Corollary~\ref{cor-Roman} gives the following.

\begin{cor}%%%%%%%%%%%%%%%%%%%%%%%%%%%%%%%%%%%%%%%%%%%%%%%%%%%%%%%%%%%%%%%%%%%%%%%%%%%%%%%%%%%%%%%%%%%%%%%%%%%%%%%%%%%%%
\label{cor-diff}
For $p\in (0,1)$, $\partial (G(n,p))\in \left\{n-2\lfloor a_{p}(n) \rfloor _{2}-i:1\leq i\leq 3\right\}$ with probability that tend to $1$ as $n\rightarrow \infty $.
\end{cor}
%%%%%%%%%%%%%%%%%%%%%%%%%%%%%%%%%%%%%%%%%%%%%%%%%%%%%%%%%%%%%%%%%%%%%%%%%%%%%%%%%%%%%%%%%%%%%%%%%%%%%%%%%%%%%%%%%%%%%%%%

%%%%%%%%%%%%%%%%%%%%%%%%%%%%%%%%%%%%%%%%%%%%%%%%%%%%%%%%%%%%%%%%%%%%%%%%%%%%%%%%%%%%%%%%%%%%%%%%%%%%%%%%%%%%%%%%%%%%%%%%
%%%%%%%%%%%%%%%%%%%%%%%%%%%%%%%%%%%%%%%%%%%%%%%%%%%%%%%%%%%%%%%%%%%%%%%%%%%%%%%%%%%%%%%%%%%%%%%%%%%%%%%%%%%%%%%%%%%%%%%%
%%%%%%%%%%%%%%%%%%%%%%%%%%%%%%%%%%%%%%%%%%%%%%%%%%%%%%%%%%%%%%%%%%%%%%%%%%%%%%%%%%%%%%%%%%%%%%%%%%%%%%%%%%%%%%%%%%%%%%%%
\section{Lemmas}\label{sec2}
%%%%%%%%%%%%%%%%%%%%%%%%%%%%%%%%%%%%%%%%%%%%%%%%%%%%%%%%%%%%%%%%%%%%%%%%%%%%%%%%%%%%%%%%%%%%%%%%%%%%%%%%%%%%%%%%%%%%%%%%
%%%%%%%%%%%%%%%%%%%%%%%%%%%%%%%%%%%%%%%%%%%%%%%%%%%%%%%%%%%%%%%%%%%%%%%%%%%%%%%%%%%%%%%%%%%%%%%%%%%%%%%%%%%%%%%%%%%%%%%%
%%%%%%%%%%%%%%%%%%%%%%%%%%%%%%%%%%%%%%%%%%%%%%%%%%%%%%%%%%%%%%%%%%%%%%%%%%%%%%%%%%%%%%%%%%%%%%%%%%%%%%%%%%%%%%%%%%%%%%%%

In this section, we prepare some lemmas which will be used in our argument.
We start with two fundamental lemmas related to the $c$-self domination concept.

\begin{lem}%%%%%%%%%%%%%%%%%%%%%%%%%%%%%%%%%%%%%%%%%%%%%%%%%%%%%%%%%%%%%%%%%%%%%%%%%%%%%%%%%%%%%%%%%%%%%%%%%%%%%%%%%%%%%
\label{lem2.1}
Let $a\in \mathbb{R}^{+}$ and $c\in (0,1)$, and let $G$ be a graph of order at least $a$.
Then $\gamma ^{c}(G)\leq a$ if and only if there exists a $c$-SDF $f:V(G)\rightarrow \{0,1,c\}$ of $G$ such that $a-1<w(f)\leq a$.
\end{lem}
%%%%%%%%%%%%%%%%%%%%%%%%%%%%%%%%%%%%%%%%%%%%%%%%%%%%%%%%%%%%%%%%%%%%%%%%%%%%%%%%%%%%%%%%%%%%%%%%%%%%%%%%%%%%%%%%%%%%%%%%
\proof
The ``if'' part is trivial.
Thus it suffices to prove the ``only if'' part.
Suppose that $\gamma ^{c}(G)\leq a$.
Then by Proposition~\ref{prop1.1}, there exists a $c$-SDF $f$ of $G$ such that $w(f)\leq a$ and $f(u)\in \{0,1,c\}$ for all $u\in V(G)$.
Choose $f$ so that $w(f)$ is as large as possible.
If $w(f)=|V(G)|$, then $w(f)=a$ because $w(f)\leq a\leq |V(G)|=w(f)$, as desired.
Thus we may assume that $w(f)<|V(G)|$.
Since $c\in (0,1)$, there exists a vertex $u_{0}\in V(G)$ such that $f(u_{0})\in \{0,c\}$.
Then the function $g:V(G)\rightarrow \{0,1,c\}$ with
$$
g(u)=
\begin{cases}
1 & (u=u_{0})\\
f(u) & (u\neq u_{0}).
\end{cases}
$$
is a $c$-SDF of $G$ and $w(g)>w(f)$.
This together with the maximality of $w(f)$ implies that $a<w(g)\leq w(f)+1$, and so $a-1<w(f)\leq a$.
\qed

\begin{lem}%%%%%%%%%%%%%%%%%%%%%%%%%%%%%%%%%%%%%%%%%%%%%%%%%%%%%%%%%%%%%%%%%%%%%%%%%%%%%%%%%%%%%%%%%%%%%%%%%%%%%%%%%%%%%
\label{lem2.2}
Let $s$ and $t$ be integers with $2\leq s\leq t-1$.
Let $G$ be a graph, and suppose that $\gamma ^{\frac{s}{t}}(G)$ is a non-integer and $\gamma ^{\frac{s}{t}}(G)\leq \lfloor \gamma ^{\frac{s}{t}}(G)\rfloor +\frac{s-1}{t}$.
Then $\gamma ^{\frac{1}{t}}(G)<\lfloor \gamma ^{\frac{s}{t}}(G)\rfloor $.
\end{lem}
%%%%%%%%%%%%%%%%%%%%%%%%%%%%%%%%%%%%%%%%%%%%%%%%%%%%%%%%%%%%%%%%%%%%%%%%%%%%%%%%%%%%%%%%%%%%%%%%%%%%%%%%%%%%%%%%%%%%%%%%
\proof
Let $f:V(G)\rightarrow \{0,1,\frac{s}{t}\}$ be an $\frac{s}{t}$-SDF of $G$ with $w(f)=\gamma ^{\frac{s}{t}}(G)$, and let $U=\{u\in V(G):f(u)=\frac{s}{t}\}$.
Since $\gamma ^{\frac{s}{t}}(G)$ is a non-integer, we have $U\neq \emptyset $.
If $|U|=1$, then $\gamma ^{\frac{s}{t}}(G)=\lfloor \gamma ^{\frac{s}{t}}(G) \rfloor +\frac{s}{t}$, which contradicts the second assumption of the lemma.
Thus $|U|\geq 2$.

Let $g:V(G)\rightarrow \{0,1,\frac{1}{t}\}$ be the function with
$$
g(u)=
\begin{cases}
\frac{1}{t} & (u\in U)\\
f(u) & (u\notin U).
\end{cases}
$$
Then $g$ is a $\frac{1}{t}$-SDF of $G$, and hence
$$
\gamma ^{\frac{1}{t}}(G)\leq w(g) = w(f)-\frac{|U|(s-1)}{t}\leq \gamma ^{\frac{s}{t}}(G)-\frac{2(s-1)}{t}\leq \lfloor \gamma ^{\frac{s}{t}}(G)\rfloor -\frac{s-1}{t}< \lfloor \gamma ^{\frac{s}{t}}(G)\rfloor,
$$
as desired.
\qed

The following lemmas are well-known (or proved by easy argument) in mathematics.

\begin{lem}[Stirling's formula]%%%%%%%%%%%%%%%%%%%%%%%%%%%%%%%%%%%%%%%%%%%%%%%%%%%%%%%%%%%%%%%%%%%%%%%%%%%%%%%%%%%%%%%%%
\label{lem2.3}
For $n\in \mathbb{Z}^{+}$, $n!\geq \sqrt{2\pi n}\left(\frac{n}{e}\right)^{n}$.
\end{lem}
%%%%%%%%%%%%%%%%%%%%%%%%%%%%%%%%%%%%%%%%%%%%%%%%%%%%%%%%%%%%%%%%%%%%%%%%%%%%%%%%%%%%%%%%%%%%%%%%%%%%%%%%%%%%%%%%%%%%%%%%

\begin{lem}%%%%%%%%%%%%%%%%%%%%%%%%%%%%%%%%%%%%%%%%%%%%%%%%%%%%%%%%%%%%%%%%%%%%%%%%%%%%%%%%%%%%%%%%%%%%%%%%%%%%%%%%%%%%%
\label{lem2.4}
For $x\geq 0$, $1-x\leq e^{-x}$.
\end{lem}
%%%%%%%%%%%%%%%%%%%%%%%%%%%%%%%%%%%%%%%%%%%%%%%%%%%%%%%%%%%%%%%%%%%%%%%%%%%%%%%%%%%%%%%%%%%%%%%%%%%%%%%%%%%%%%%%%%%%%%%%

%%%%%%%%%%%%%%%%%%%%%%%%%%%%%%%%%%%%%%%%%%%%%%%%%%%%%%%%%%%%%%%%%%%%%%%%%%%%%%%%%%%%%%%%%%%%%%%%%%%%%%%%%%%%%%%%%%%%%%%%
%%%%%%%%%%%%%%%%%%%%%%%%%%%%%%%%%%%%%%%%%%%%%%%%%%%%%%%%%%%%%%%%%%%%%%%%%%%%%%%%%%%%%%%%%%%%%%%%%%%%%%%%%%%%%%%%%%%%%%%%
%%%%%%%%%%%%%%%%%%%%%%%%%%%%%%%%%%%%%%%%%%%%%%%%%%%%%%%%%%%%%%%%%%%%%%%%%%%%%%%%%%%%%%%%%%%%%%%%%%%%%%%%%%%%%%%%%%%%%%%%
\section{A crude estimation}\label{sec3}
%%%%%%%%%%%%%%%%%%%%%%%%%%%%%%%%%%%%%%%%%%%%%%%%%%%%%%%%%%%%%%%%%%%%%%%%%%%%%%%%%%%%%%%%%%%%%%%%%%%%%%%%%%%%%%%%%%%%%%%%
%%%%%%%%%%%%%%%%%%%%%%%%%%%%%%%%%%%%%%%%%%%%%%%%%%%%%%%%%%%%%%%%%%%%%%%%%%%%%%%%%%%%%%%%%%%%%%%%%%%%%%%%%%%%%%%%%%%%%%%%
%%%%%%%%%%%%%%%%%%%%%%%%%%%%%%%%%%%%%%%%%%%%%%%%%%%%%%%%%%%%%%%%%%%%%%%%%%%%%%%%%%%%%%%%%%%%%%%%%%%%%%%%%%%%%%%%%%%%%%%%

In this section, we prove the following theorem which is weaker than Theorem~\ref{thm-main}

\begin{thm}%%%%%%%%%%%%%%%%%%%%%%%%%%%%%%%%%%%%%%%%%%%%%%%%%%%%%%%%%%%%%%%%%%%%%%%%%%%%%%%%%%%%%%%%%%%%%%%%%%%%%%%%%%%%%
\label{thm3.0}
Let $s$ and $t$ be integers with $1\leq s\leq t-1$.
Then for $p\in (0,1)$,
$$
\gamma ^{\frac{t}{s}}(G(n,p))\in \left[\lfloor a_{p}(n) \rfloor _{t}+\frac{1}{t},~b_{p,t}(n)\right]
$$
with probability that tend to $1$ as $n\rightarrow \infty $.
\end{thm}
%%%%%%%%%%%%%%%%%%%%%%%%%%%%%%%%%%%%%%%%%%%%%%%%%%%%%%%%%%%%%%%%%%%%%%%%%%%%%%%%%%%%%%%%%%%%%%%%%%%%%%%%%%%%%%%%%%%%%%%%

In \cite{WG}, Wieland and Godbole implicitly proved the following lemma.

\begin{lem}[Wieland and Godbole~\cite{WG}]%%%%%%%%%%%%%%%%%%%%%%%%%%%%%%%%%%%%%%%%%%%%%%%%%%%%%%%%%%%%%%%%%%%%%%%%%%%%%%
\label{ThmC}
Let $\varepsilon \in \mathbb{R}^{+}$.
Then for $p\in (0,1)$, $\gamma (G(n,p))\leq \lceil a_{p}(n)+\varepsilon \rceil $ with probability that tend to $1$ as $n\rightarrow \infty $.
\end{lem}
%%%%%%%%%%%%%%%%%%%%%%%%%%%%%%%%%%%%%%%%%%%%%%%%%%%%%%%%%%%%%%%%%%%%%%%%%%%%%%%%%%%%%%%%%%%%%%%%%%%%%%%%%%%%%%%%%%%%%%%%

\begin{lem}%%%%%%%%%%%%%%%%%%%%%%%%%%%%%%%%%%%%%%%%%%%%%%%%%%%%%%%%%%%%%%%%%%%%%%%%%%%%%%%%%%%%%%%%%%%%%%%%%%%%%%%%%%%%%
\label{lem3.00}
For $p\in (0,1)$, $t\in \mathbb{Z}^{+}$ and $n\in \mathbb{Z}^{+}\setminus \{1\}$, we have $\lceil a_{p}(n)+\frac{1}{2t} \rceil \leq b_{p,t}(n)$.
\end{lem}
%%%%%%%%%%%%%%%%%%%%%%%%%%%%%%%%%%%%%%%%%%%%%%%%%%%%%%%%%%%%%%%%%%%%%%%%%%%%%%%%%%%%%%%%%%%%%%%%%%%%%%%%%%%%%%%%%%%%%%%%
\proof
There exist non-negative integers $m_{1}$ and $m_{2}$ such that $m_{1}+\frac{m_{2}}{t}\leq a_{p}(n)<m_{1}+\frac{m_{2}+1}{t}$ and $0\leq m_{2}\leq t-1$.

Suppose that $m_{2}=t-1$.
Since $\lfloor a_{p}(n) \rfloor _{t}+\frac{1}{t}=m_{1}+\frac{t-1}{t}+\frac{1}{t}=m_{1}+1~(\in \mathbb{Z}^{+})$, we have $b_{p,t}(n)=\lfloor \lfloor a_{p}(n) \rfloor _{t}+\frac{1}{t} \rfloor +1=m_{1}+2$.
On the other hand, $a_{p}(n)+\frac{1}{2t} < m_{1}+1+\frac{1}{2t}$, and so $\lceil a_{p}(n)+\frac{1}{2t} \rceil \leq m_{1}+2=b_{p,t}(n)$, as desired.
Thus we may assume that $0\leq m_{2}\leq t-2$.

Since $\lfloor a_{p}(n) \rfloor _{t}+\frac{1}{t}=m_{1}+\frac{m_{2}+1}{t}\leq m_{1}+\frac{t-1}{t}$, we have $b_{p,t}(n)=\lfloor \lfloor a_{p}(n) \rfloor _{t}+\frac{1}{t} \rfloor +1=m_{1}+1$.
On the other hand, $a_{p}(n)+\frac{1}{2t} < m_{1}+\frac{t-1}{t}+\frac{1}{2t}=m_{1}+\frac{2t-1}{2t}<m_{1}+1$, and so $\lceil a_{p}(n)+\frac{1}{2t} \rceil \leq m_{1}+1=b_{p,t}(n)$, as desired.
\qed

\medbreak\noindent\textit{Proof of Theorem~\ref{thm3.0}.}\quad
Note that $\gamma ^{\frac{s}{t}}(G)\leq \gamma ^{1}(G)=\gamma (G)$ for all graphs $G$.
Hence by Lemma~\ref{ThmC} with $\varepsilon =\frac{1}{2t}$ and Lemma~\ref{lem3.00},
\begin{align*}
{\rm Pr}(\gamma ^{\frac{s}{t}}(G(n,p))\leq b_{p,t}(n)) &\geq {\rm Pr}(\gamma (G(n,p))\leq b_{p,t}(n))\\
&\geq {\rm Pr}\left(\gamma (G(n,p))\leq \left\lceil a_{p}(n)+\frac{1}{2t} \right\rceil \right)\\
&\rightarrow 1~~(n\rightarrow \infty ).
\end{align*}
Consequently, we obtain the upper bound of the theorem.

We next prove the lower bound of the theorem.
Let $\mathbb{M}=\{m_{1}+\frac{m_{2}}{t}:m_{1},m_{2}\in \mathbb{Z}^{+}\cup \{0\}\}$, and for $a\in \mathbb{R}^{+}$, let $\MM(a)=\{(m_{1},m_{2}):m_{1}+\frac{m_{2}}{t}=a\}$.
Then $\MM(a)\neq \emptyset $ if and only if $a\in \mathbb{M}$.
Furthermore, we note that $\lfloor a_{p}(n) \rfloor _{t}+\frac{1}{t}$ is the smallest number in $\mathbb{M}$ more than $a_{p}(n)$.
Since $\gamma ^{\frac{s}{t}}(G)\geq \gamma ^{\frac{1}{t}}(G)$ for all graphs $G$, it suffices to show that $\gamma ^{\frac{1}{t}}(G(n,p))>a_{p}(n)$ with probability that tend to $1$ as $n\rightarrow \infty $.

%We consider a random variable concerning dominating functions of random graphs $G(n,p)$.
For $m_{1},m_{2}\in \mathbb{Z}^{+}\cup \{0\}$, let $X_{m_{1},m_{2}}$ be the random variable counting the number of $\frac{1}{t}$-SDFs $f:[n]\rightarrow \{0,1,\frac{1}{t}\}$ of $G(n,p)$ with $|\{u\in [n]:f(u)=1\}|=m_{1}$ and $|\{u\in [n]:f(u)=\frac{1}{t}\}|=m_{2}$.
For $a\in \mathbb{M}$, let $X_{a}=\sum _{(m_{1},m_{2})\in \MM(a)}X_{m_{1},m_{2}}$.

For a graph $G$, an ordered pair $(S_{1},S_{2})$ of subsets of $V(G)$ with $S_{1}\cap S_{2}=\emptyset $ is called a {\it $\frac{1}{t}$-self dominating pair} of $G$ if the function $f:V(G)\rightarrow \{0,1,\frac{1}{t}\}$ with
$$
f(u)=
\begin{cases}
0 & (u\in V(G)\setminus (S_{1}\cup S_{2}))\\
1 & (u\in S_{1})\\
\frac{1}{t} & (u\in S_{2})
\end{cases}
$$
is a $\frac{1}{t}$-SDF of $G$.
Let $\SS _{m_{1},m_{2}}=\left\{(S_{1},S_{2})\in \binom{[n]}{m_{1}}\times \binom{[n]}{m_{2}}:S_{1}\cap S_{2}=\emptyset \right\}$, and for $(S_{1},S_{2})\in \SS_{m_{1},m_{2}}$, let $I_{S_{1},S_{2}}$ be the random variable satisfying
$$
I_{S_{1},S_{2}}=
\begin{cases}
1 & (\mbox{$(S_{1},S_{2})$ is a $\frac{1}{t}$-self dominating pair of $G(n,p)$})\\
0 & (\mbox{otherwise}).
\end{cases}
$$
Note that $X_{m_{1},m_{2}}=\sum _{(S_{1},S_{2})\in \SS_{m_{1},m_{2}}}I_{S_{1},S_{2}}$.
The following claim plays a key role in our argument.

\begin{claim}%%%%%%%%%%%%%%%%%%%%%%%%%%%%%%%%%%%%%%%%%%%%%%%%%%%%%%%%%%%%%%%%%%%%%%%%%%%%%%%%%%%%%%%%%%%%%%%%%%%%%%%%%%%
\label{cl3.1}
For non-negative integers $m_{1}$ and $m_{2}$, $E(X_{m_{1},m_{2}})=\frac{n!}{(n-m_{1}-m_{2})!~m_{1}!~m_{2}!}(1-(1-p)^{m_{1}})^{n-m_{1}-m_{2}}$.
\end{claim}
%%%%%%%%%%%%%%%%%%%%%%%%%%%%%%%%%%%%%%%%%%%%%%%%%%%%%%%%%%%%%%%%%%%%%%%%%%%%%%%%%%%%%%%%%%%%%%%%%%%%%%%%%%%%%%%%%%%%%%%%
\proof
For $(S_{1},S_{2})\in \SS_{m_{1},m_{2}}$, since ${\rm Pr}(N_{G}(u)\cap S_{1}\neq \emptyset )=1-(1-p)^{m_{1}}$ for each $u\in [n]\setminus (S_{1}\cup S_{2})$,
$$
{\rm Pr}(I_{S_{1},S_{2}}=1)=\prod _{u\in [n]\setminus (S_{1}\cup S_{2})}{\rm Pr}(N_{G}(u)\cap S_{1}\neq \emptyset )=(1-(1-p)^{m_{1}})^{n-m_{1}-m_{2}}.
$$
Since $X_{m_{1},m_{2}}=\sum _{(S_{1},S_{2})\in \SS_{m_{1},m_{2}}}I_{S_{1},S_{2}}$, it follows that
\begin{align*}
E(X_{m_{1},m_{2}}) &= \sum _{(S_{1},S_{2})\in \SS_{m_{1},m_{2}}}E(I_{S_{1},S_{2}})\\
&= \sum _{(S_{1},S_{2})\in \SS_{m_{1},m_{2}}}{\rm Pr}(I_{S_{1},S_{2}}=1)\\
&= \binom{n}{m_{1}}\binom{n-m_{1}}{m_{2}}(1-(1-p)^{m_{1}})^{n-m_{1}-m_{2}},
\end{align*}
as desired.
\qed

Since $\frac{1}{1-p}>1$, the value $h_{0}=\min\{h\in \mathbb{Z}^{+}:t-\frac{1}{(1-p)^{a}}<0$ for all $a\geq h\}$ is a well-defined constant (depending on $p$ and $t$ only).
In the rest of this proof, we consider $G(n,p)$ for sufficiently large $n$.
Thus, for example, we may assume that $L(L(n))>0$, $n>ta_{p}(n)$, $a_{p}(n)>h_{0}$, etc.
For $x\in \mathbb{R}^{+}$, let $L(x)=\log _{1/(1-p)}x$.
Note that $a_{p}(n)=\log _{1/(1-p)}\frac{n}{\log _{1/(1-p)}n\ln n}=L(\frac{n}{L(n)\ln n})$.

\begin{claim}%%%%%%%%%%%%%%%%%%%%%%%%%%%%%%%%%%%%%%%%%%%%%%%%%%%%%%%%%%%%%%%%%%%%%%%%%%%%%%%%%%%%%%%%%%%%%%%%%%%%%%%%%%%
\label{cl-low-1}
Let $m_{1}$ and $m_{2}$ be non-negative integers with $a_{p}(n)-1<m_{1}+\frac{m_{2}}{t}\leq a_{p}(n)$.
Then the following hold.
\begin{enumerate}
\item[{\upshape(i)}]
We have $E(X_{m_{1},m_{2}})<\exp\left[(m_{1}+m_{2})(\ln n+2)-\frac{L(n)\ln n}{(1-p)^{a_{p}(n)-m_{1}}}\right]$.
\item[{\upshape(ii)}]
If $0\leq m_{1}\leq a_{p}(n)-h_{0}$, then $E(X_{m_{1},m_{2}})<\exp[t(2L(n)-L(L(n)\ln n)\ln n)]$.
\end{enumerate}
\end{claim}
%%%%%%%%%%%%%%%%%%%%%%%%%%%%%%%%%%%%%%%%%%%%%%%%%%%%%%%%%%%%%%%%%%%%%%%%%%%%%%%%%%%%%%%%%%%%%%%%%%%%%%%%%%%%%%%%%%%%%%%%
\proof
\begin{enumerate}
\item[{\upshape(i)}]
By Lemma~\ref{lem2.3}, if $m_{1}\geq 1$ and $m_{2}\geq 1$, then
$$
\frac{n!}{(n-m_{1}-m_{2})!~m_{1}!~m_{2}!} \leq n^{m_{1}+m_{2}}\cdot \frac{1}{\sqrt{2\pi m_{1}}\left(\frac{m_{1}}{e}\right)^{m_{1}}}\cdot \frac{1}{\sqrt{2\pi m_{2}}\left(\frac{m_{2}}{e}\right)^{m_{2}}} < (en)^{m_{1}+m_{2}};
$$
if $m_{i}=0$ for some $i\in \{1,2\}$, then $m_{3-i}\geq 1$, and hence
$$
\frac{n!}{(n-m_{1}-m_{2})!~m_{1}!~m_{2}!} \leq n^{m_{3-i}}\cdot \frac{1}{\sqrt{2\pi m_{3-i}}\left(\frac{m_{3-i}}{e}\right)^{m_{3-i}}}<(en)^{m_{3-i}}=(en)^{m_{1}+m_{2}}.
$$
In either case,
\begin{align}
\frac{n!}{(n-m_{1}-m_{2})!~m_{1}!~m_{2}!} < (en)^{m_{1}+m_{2}}.\label{eq-low-1}
\end{align}
Furthermore, we have
\begin{align}
n(1-p)^{m_{1}} = \frac{n(1-p)^{L(\frac{n}{L(n)\ln n})}}{(1-p)^{a_{p}(n)-m_{1}}} = \frac{n\cdot \frac{L(n)\ln n}{n}}{(1-p)^{a_{p}(n)-m_{1}}} = \frac{L(n)\ln n}{(1-p)^{a_{p}(n)-m_{1}}}.\label{eq-low-2}
\end{align}
By Claim~\ref{cl3.1}, Lemma~\ref{lem2.4}, (\ref{eq-low-1}) and (\ref{eq-low-2}),
\begin{align*}
E(X_{m_{1},m_{2}}) &= \frac{n!}{(n-m_{1}-m_{2})!~m_{1}!~m_{2}!}(1-(1-p)^{m_{1}})^{n-m_{1}-m_{2}}\\
&< (en)^{m_{1}+m_{2}}\left(e^{-(1-p)^{m_{1}}}\right)^{n-m_{1}-m_{2}}\\
&= \exp[(m_{1}+m_{2})+(m_{1}+m_{2})\ln n-n(1-p)^{m_{1}}+(m_{1}+m_{2})(1-p)^{m_{1}}]\\
&\leq \exp\left[2(m_{1}+m_{2})+(m_{1}+m_{2})\ln n-\frac{L(n)\ln n}{(1-p)^{a_{p}(n)-m_{1}}}\right].
\end{align*}

\item[{\upshape(ii)}]
By the definition of $m_{1}$ and $m_{2}$, we have
\begin{align}
m_{1}+m_{2}\leq t\left(m_{1}+\frac{m_{2}}{t}\right)\leq ta_{p}(n)=t(L(n)-L(L(n)\ln n)).\label{eq-low-3}
\end{align}
Since $a_{p}(n)-m_{1}\geq h_{0}$, it follows from the definition of $h_{0}$ that $(t-\frac{1}{(1-p)^{a_{p}(n)-m_{1}}})L(n)\ln n<0$.
This together with (i) and (\ref{eq-low-3}) implies that
\begin{align*}
E(X_{m_{1},m_{2}}) &\leq \exp\left[(m_{1}+m_{2})(\ln n+2)-\frac{L(n)\ln n}{(1-p)^{a_{p}(n)-m_{1}}}\right]\\
&\leq \exp\left[t(L(n)-L(L(n)\ln n))(\ln n+2)-\frac{L(n)\ln n}{(1-p)^{a_{p}(n)-m_{1}}}\right]\\
&= \exp\left[\left(t-\frac{1}{(1-p)^{a_{p}(n)-m_{1}}} \right)L(n) \ln n+t(2L(n)-L(L(n)\ln n)\ln n-2L(L(n)\ln n))\right]\\
&< \exp[t(2L(n)-L(L(n)\ln n)\ln n)],
\end{align*}
as desired.
\qed
\end{enumerate}

\begin{claim}%%%%%%%%%%%%%%%%%%%%%%%%%%%%%%%%%%%%%%%%%%%%%%%%%%%%%%%%%%%%%%%%%%%%%%%%%%%%%%%%%%%%%%%%%%%%%%%%%%%%%%%%%%%
\label{cl-low-2}
Let $a\in \mathbb{M}$ be a number with $a_{p}(n)-1<a\leq a_{p}(n)$.
Then $E(X_{a})\rightarrow 0$ if $n\rightarrow \infty $.
\end{claim}
%%%%%%%%%%%%%%%%%%%%%%%%%%%%%%%%%%%%%%%%%%%%%%%%%%%%%%%%%%%%%%%%%%%%%%%%%%%%%%%%%%%%%%%%%%%%%%%%%%%%%%%%%%%%%%%%%%%%%%%%
\proof
By the definition of $X_{a}$,
\begin{align*}
E(X_{a}) = E\left(\sum _{(m_{1},m_{2})\in \MM(a)}X_{m_{1},m_{2}}\right)= \sum _{\substack{(m_{1},m_{2})\in \MM(a)\\0\leq m_{1}\leq a_{p}(n)-h_{0}}}E(X_{m_{1},m_{2}})+\sum _{\substack{(m_{1},m_{2})\in \MM(a)\\a_{p}(n)-h_{0}<m_{1}\leq a}}E(X_{m_{1},m_{2}}).
%&= \sum _{(m_{1},m_{2})\in \MM_{\frac{s}{t}}(a)}E(X_{m_{1},m_{2}})\\
%&\leq \sum _{\substack{(m_{1},m_{2})\in (\mathbb{Z}^{+}\cup \{0\})^{2}\\m_{1}\leq a}}E(X_{m_{1},m_{2}})\\
\end{align*}
Note that the number of $m_{1}\in \mathbb{Z}^{+}$ satisfying $a_{p}(n)-h_{0}<m_{1}\leq a$ is at most $h_{0}$ because $a\leq a_{p}(n)$.
Hence $\sum _{\substack{(m_{1},m_{2})\in \MM(a)\\a_{p}(n)-h_{0}<m_{1}\leq a}}E(X_{m_{1},m_{2}})$ is a sum having constant terms.
Thus it suffices to prove the following:
\begin{enumerate}
\item[{\bf (A1)}]
$\sum _{\substack{(m_{1},m_{2})\in \MM(a)\\0\leq m_{1}\leq a_{p}(n)-h_{0}}}E(X_{m_{1},m_{2}})\rightarrow 0~~(n\rightarrow \infty )$, and
\item[{\bf (A2)}]
for each $(m_{1},m_{2})\in \MM(a)$, if $a_{p}(n)-h_{0}<m_{1}\leq a$, then $E(X_{m_{1},m_{2}})\rightarrow 0~~(n\rightarrow \infty )$.
\end{enumerate}

By Claim~\ref{cl-low-1}(ii),
\begin{align*}
\sum _{\substack{(m_{1},m_{2})\in \MM(a)\\0\leq m_{1}\leq a_{p}(n)-h_{0}}}E(X_{m_{1},m_{2}}) &< (a_{p}(n)-h_{0}+1)\exp[t(2L(n)-L(L(n)\ln n)\ln n)]\\
&\leq a_{p}(n)\exp[t(2L(n)-L(L(n)\ln n)\ln n)]\\
&= \exp[\ln a_{p}(n)+t(2L(n)-L(L(n)\ln n)\ln n)]\\
&< \exp[\ln L(n)+t(2L(n)-L(L(n)\ln n)\ln n)]\\
&\rightarrow 0~~(n\rightarrow \infty ),
\end{align*}
which proves (A1).

We next assume that $(m_{1},m_{2})\in \MM(a)$ satisfies $a_{p}(n)-h_{0}<m_{1}\leq a$ and prove (A2).
We have
$$
m_{1}+m_{2} = t\left(m_{1}+\frac{1}{t}m_{2}\right)-(t-1)m_{1}< ta_{p}(n)-(t-1)(a_{p}(n)-h_{0})= a_{p}(n)+(t-1)h_{0}.
$$
Note that $\alpha :=(t-1)h_{0}$ is a constant depending on $p$ and $t$ only.
Hence it follows from Claim~\ref{cl-low-1}(i) that
\begin{align*}
E(X_{m_{1},m_{2}}) &\leq \exp\left[(m_{1}+m_{2})(\ln n+2)-\frac{L(n)\ln n}{(1-p)^{a_{p}(n)-m_{1}}}\right]\\
&< \exp[(a_{p}(n)+\alpha )(\ln n+2)-L(n)\ln n]\\
&= \exp[-L(L(n)\ln n)\ln n+2L(n)-2L(L(n)\ln n)+\alpha \ln n+2\alpha ]\\
&\rightarrow 0~~(n\rightarrow \infty ),
\end{align*}
which proves (A2).
\qed

Let $A_{n}=\{a\in \mathbb{M}:a_{p}(n)-1<a\leq a_{p}(n)\}$.
Then $|A_{n}|\leq t$.
In particular, $\sum _{a\in A_{n}}E(X_{a})$ is a sum having constant terms.
Consequently, it follows from Lemma~\ref{lem2.1} and Claim~\ref{cl-low-2} that
$$
{\rm Pr}(\gamma ^{\frac{1}{t}}(G(n,p))\leq a_{p}(n))\leq \sum _{a\in A_{n}}{\rm Pr}(X_{a}\geq 1)\leq \sum _{a\in A_{n}}E(X_{a})\rightarrow 0~~(n\rightarrow \infty ),
$$
and so ${\rm Pr}(\gamma ^{\frac{1}{t}}(G(n,p))>a_{p}(n))\rightarrow 1~~(n\rightarrow \infty )$.

This completes the proof of Theorem~\ref{thm3.0}.
\qed

%%%%%%%%%%%%%%%%%%%%%%%%%%%%%%%%%%%%%%%%%%%%%%%%%%%%%%%%%%%%%%%%%%%%%%%%%%%%%%%%%%%%%%%%%%%%%%%%%%%%%%%%%%%%%%%%%%%%%%%%
%%%%%%%%%%%%%%%%%%%%%%%%%%%%%%%%%%%%%%%%%%%%%%%%%%%%%%%%%%%%%%%%%%%%%%%%%%%%%%%%%%%%%%%%%%%%%%%%%%%%%%%%%%%%%%%%%%%%%%%%
%%%%%%%%%%%%%%%%%%%%%%%%%%%%%%%%%%%%%%%%%%%%%%%%%%%%%%%%%%%%%%%%%%%%%%%%%%%%%%%%%%%%%%%%%%%%%%%%%%%%%%%%%%%%%%%%%%%%%%%%
\section{Graph-theoretical refinement of Theorem~\ref{thm3.0}}\label{sec4}
%%%%%%%%%%%%%%%%%%%%%%%%%%%%%%%%%%%%%%%%%%%%%%%%%%%%%%%%%%%%%%%%%%%%%%%%%%%%%%%%%%%%%%%%%%%%%%%%%%%%%%%%%%%%%%%%%%%%%%%%
%%%%%%%%%%%%%%%%%%%%%%%%%%%%%%%%%%%%%%%%%%%%%%%%%%%%%%%%%%%%%%%%%%%%%%%%%%%%%%%%%%%%%%%%%%%%%%%%%%%%%%%%%%%%%%%%%%%%%%%%
%%%%%%%%%%%%%%%%%%%%%%%%%%%%%%%%%%%%%%%%%%%%%%%%%%%%%%%%%%%%%%%%%%%%%%%%%%%%%%%%%%%%%%%%%%%%%%%%%%%%%%%%%%%%%%%%%%%%%%%%

In this section, we complete the proof of Theorem~\ref{thm-main}.
Let $s$, $t$ and $p$ be numbers as in Theorem~\ref{thm-main}.
Let $\varepsilon \in \mathbb{R}^{+}$.
Then by Theorem~\ref{thm3.0}, there exists $N_{0}\in \mathbb{Z}^{+}$ such that for every integer $n\geq N_{0}$,
$$
{\rm Pr}\left(\gamma ^{\frac{1}{t}}(G(n,p))<\lfloor a_{p}(n) \rfloor _{t}+\frac{1}{t}\right)<\frac{\varepsilon }{2(s-1)}~~~\mbox{and}~~~{\rm Pr}\left(\gamma ^{\frac{t}{s}}(G(n,p))\notin \left[\lfloor a_{p}(n) \rfloor _{t}+\frac{1}{t},~b_{p,t}(n)\right]\right)<\frac{\varepsilon }{2}.
$$

Fix an integer $n\geq N_{0}$, and let $i$ be an integer with $t-s+1\leq i\leq t-1$.
Since $b_{p,t}(n)$ is an integer, $b_{p,t}(n)-\frac{i}{t}$ is a non-integer.
Furthermore, if a graph $G$ satisfies $\gamma ^{\frac{s}{t}}(G)=b_{p,t}(n)-\frac{i}{t}$, then
$$
\lfloor \gamma ^{\frac{s}{t}}(G) \rfloor = \left\lfloor b_{p,t}(n)-\frac{i}{t} \right\rfloor =b_{p,t}(n)-1,
$$
and hence
$$
\gamma ^{\frac{s}{t}}(G) = b_{p,t}(n)-\frac{i}{t}= \lfloor \gamma ^{\frac{s}{t}}(G) \rfloor +1-\frac{i}{t}
\leq \lfloor \gamma ^{\frac{s}{t}}(G) \rfloor +1-\frac{t-s+1}{t}=\lfloor \gamma ^{\frac{s}{t}}(G) \rfloor +\frac{s-1}{t}.
$$
This together with Lemma~\ref{lem2.2} implies that if $\gamma ^{\frac{s}{t}}(G)=b_{p,t}(n)-\frac{i}{t}$, then $\gamma ^{\frac{1}{t}}(G)<\lfloor \gamma ^{\frac{s}{t}}(G)\rfloor =b_{p,t}(n)-1$.
Hence we have ${\rm Pr}(\gamma ^{\frac{1}{t}}(G(n,p))<b_{p,t}(n)-1)\geq {\rm Pr}(\gamma ^{\frac{s}{t}}(G(n,p))=b_{p,t}(n)-\frac{i}{t})$.
On the other hand, since $b_{p,t}(n)-1=\lfloor \lfloor a_{p}(n) \rfloor _{t}+\frac{1}{t} \rfloor \leq \lfloor a_{p}(n) \rfloor _{t}+\frac{1}{t}$,
\begin{align*}
{\rm Pr}\left(\gamma ^{\frac{s}{t}}(G(n,p))=b_{p,t}(n)-\frac{i}{t}\right) &\leq {\rm Pr}(\gamma ^{\frac{1}{t}}(G(n,p))<b_{p,t}(n)-1)\\
&\leq {\rm Pr}\left(\gamma ^{\frac{1}{t}}(G(n,p))<\lfloor a_{p}(n) \rfloor _{t}+\frac{1}{t}\right)\\
&< \frac{\varepsilon }{2(s-1)}.
\end{align*}
Consequently,
$$
{\rm Pr}\left(\gamma ^{\frac{s}{t}}(G(n,p))\in \left\{b_{p,t}(n)-\frac{i}{t}:t-s+1\leq i\leq t-1\right\}\right)<\frac{\varepsilon }{2},
$$
and so
$$
{\rm Pr}\left(\gamma ^{\frac{s}{t}}(G(n,p))\notin \left[\lfloor a_{p}(n) \rfloor _{t}+\frac{1}{t},~b_{p,t}(n)\right]\mbox{ or }\gamma ^{\frac{s}{t}}(G(n,p))\in \left\{b_{p,t}(n)-\frac{i}{t}:t-s+1\leq i\leq t-1\right\}\right)<\frac{\varepsilon }{2}+\frac{\varepsilon }{2}=\varepsilon .
$$
Since $\varepsilon $ is arbitrary, this completes the proof of Theorem~\ref{thm-main}.

%%%%%%%%%%%%%%%%%%%%%%%%%%%%%%%%%%%%%%%%%%%%%%%%%%%%%%%%%%%%%%%%%%%%%%%%%%%%%%%%%%%%%%%%%%%%%%%%%%%%%%%%%%%%%%%%%%%%%%%%
%%%%%%%%%%%%%%%%%%%%%%%%%%%%%%%%%%%%%%%%%%%%%%%%%%%%%%%%%%%%%%%%%%%%%%%%%%%%%%%%%%%%%%%%%%%%%%%%%%%%%%%%%%%%%%%%%%%%%%%%
%%%%%%%%%%%%%%%%%%%%%%%%%%%%%%%%%%%%%%%%%%%%%%%%%%%%%%%%%%%%%%%%%%%%%%%%%%%%%%%%%%%%%%%%%%%%%%%%%%%%%%%%%%%%%%%%%%%%%%%%
\section*{Acknowledgment}
%%%%%%%%%%%%%%%%%%%%%%%%%%%%%%%%%%%%%%%%%%%%%%%%%%%%%%%%%%%%%%%%%%%%%%%%%%%%%%%%%%%%%%%%%%%%%%%%%%%%%%%%%%%%%%%%%%%%%%%%
%%%%%%%%%%%%%%%%%%%%%%%%%%%%%%%%%%%%%%%%%%%%%%%%%%%%%%%%%%%%%%%%%%%%%%%%%%%%%%%%%%%%%%%%%%%%%%%%%%%%%%%%%%%%%%%%%%%%%%%%
%%%%%%%%%%%%%%%%%%%%%%%%%%%%%%%%%%%%%%%%%%%%%%%%%%%%%%%%%%%%%%%%%%%%%%%%%%%%%%%%%%%%%%%%%%%%%%%%%%%%%%%%%%%%%%%%%%%%%%%%

The authors would like to thank Professor Yoshimi Egawa for his helpful comments about Section~\ref{sec4}.
This work was partially supported by JSPS KAKENHI Grant number JP18K13449 (to M.F).

\end{document}